\documentclass{article}
\usepackage{amssymb}

\usepackage{amsmath}
\usepackage{chicago}

\newtheorem{theorem}{Theorem}[section]

\makeatletter
\@addtoreset{equation}{section}
\makeatother

\newtheorem{conjecture}{Conjecture}[section]
\newtheorem{corollary}{Corollary}[section]

\newtheorem{example}[theorem]{Example}

\newtheorem{lemma}{Lemma}[section]

\newtheorem{proposition}[theorem]{Proposition}

\input{tcilatex}

\begin{document}

\title{Uncertainty of the Shapley Value}
\author{Vladislav Kargin\thanks{%
Cornerstone Research, 599 Lexington Avenue, New York, NY 10022, USA;
slava@bu.edu} \thanks{%
The views expressed within this article are the views of the author and not
the views of Cornerstone Research.}}
\maketitle

\begin{abstract}
This paper introduces a measure of uncertainty in the determination of the
Shapley value, illustrates it with examples, and studies some of its
properties. The introduced measure of uncertainty quantifies random
variations in a player's marginal contribution during the bargaining
process. The measure is symmetric with respect to exchangeable substitutions
in the players, equal to zero for dummy player, and convex in the game
argument. The measure is illustrated by several examples of abstract games
and an example from epidemiology.
\end{abstract}

\begin{quotation}
\newpage

The dread of beatings! Dread of being late!

And, greatest dread of all, the dread of games!

Sir John Betjeman, ``Summoned by Bells'', ch.7
\end{quotation}

\section{Introduction}

An important problem in epidemiology is to assign weights to the factors
that contribute to an increase in the incidence rate of a disease. For
example, a higher incidence rate of a heart disease may be attributable to
smoking, bad eating habits, or lack of physical exercise, and
epidemiologists are interested in quantifying the relative degree of
importance of these factors. There are different methods to perform the
attribution, all of which must deal with the fact that the contribution of a
factor to the incidence rate depends on the presence of the other factors.
For instance, smoking may hypothetically have a relatively small impact on
health unless it is accompanied by other factors. This difficulty is similar
to the problems that arise in cooperative game theory, and in a recent paper %
\citeANP{land_gefeller97} have proposed an attribution according to the
Shapley value, the standard solution concept from game theory. Specifically,
they have studied the influence of smoking and three types of cholesterol
level, LDL, HDL, and VLDL, on a certain heart disease, myocardial
infarction. When all four factors are absent, the incidence rate of the
disease is $77$ percent lower than the incidence rate in the total
population. The Shapley attribution for this reduction is $21$ percent for
smoking, and $41,$ $8,$ and $8$ percent for the cholesterol levels
respectively. A natural question arises about the statistical significance
of these attributions. The present paper proposes a measure of the
uncertainty of the Shapley value based on the probabilistic interpretation
of this game-theoretic concept.

The uncertainty of the Shapley value is going to be useful not only in
practical applications similar to the epidemiological example above but also
in theoretical problems. For example, following a line of research initiated
by Shapley, \citeN{aumann75} has shown that as the number of participants in
a market economy grows, the Shapley value of a game associated with this
economy converges to the competitive equilibrium allocation. This
convergence puts the concept of the competitive equilibrium on firm
game-theoretical grounds, but only if it is presumed that the Shapley value
reflects the results of the bargaining in the economy accurately. If the
Shapley value provides the outcome of the bargaining with only a limited
degree of certainty, then the natural question arises: What if the
uncertainty of the Shapley value increases with the increase in the number
of participants? Wouldn't it imply that the outcome of bargaining does not
necessarily converge to competitive equilibrium? And if the uncertainty does
decrease, then how does the rate of the decrease depend on the details of
the market organization? Answering \ these questions is important for a
better understanding of how the competitive equilibrium is achieved in
market economies.

From another perspective, the Shapley value is used not only as a device to
compute a fair allocation of resources among players but also as a tool to
evaluate the prospect of playing a game. The bargaining uncertainty would
provide an additional dimension for evaluating this prospect. Indeed, both
experimental evidence and personal introspection suggest that participation
in a bargaining situation implies undergoing the vagaries of the bargaining
process, an experience that may be a source of pleasure or displeasure
depending on the personality of the participant. A measure of bargaining
uncertainty would serve as a useful tool to investigate this aspect of a
game. Characterizing a game by both its value and uncertainty is similar to
characterizing a weapon by its power and precision, or a financial stock by
its expected return and risk. In all these situations the complexity of a
real-world object cannot be compressed into just one dimension, it needs at
least two for proper description.

So, how should the game uncertainty be measured? To answer this question,
let us first examine the different frameworks that are used to define the
Shapley value.

The shortest way to define the Shapley value is axiomatic. However, this
approach is not suitable for measuring the Shapley value uncertainty. This
is because the Shapley value is the unique solution satisfying all the
axioms, and relaxation of any of them leads to an infinity of new solutions.
In contrast, the probabilistic interpretation of the Shapley value suits the
problem of defining uncertainty very well: the Shapley value is defined as
the expectation of the marginal contribution to a random coalition, and it
is straightforward to define the uncertainty as the variance of the marginal
contribution. This is the way of defining uncertainty that the present paper
propounds.

The probabilistic interpretation of the Shapley value was sketched by
Shapley in his seminal paper (\citeyear{shapley53}). Later it was set on
firm ground by \citeN{gul89}, \citeN{evans96}, and \citeN{hart_mas-colell96}
in their work on non-cooperative foundations of the Shapley value, work that
was inspired by Rubinstein's paper on the bargaining problem (%
\citeyear{rubinstein82}). Essentially, in all these probabilistic
interpretations, the randomness of the model is the randomness of bargaining
order and coalition formation.

To understand how the randomness in the order of bargaining introduces
uncertainty about the outcome, it is useful to recall Rubinstein's model of
the bargaining game. In this game two players aim to divide a pie that has
unit value. The players are impatient and make offers in turn until some of
them agrees to the other player's offer. The player who makes the first
offer is determined by chance. Rubinstein has shown that there is a unique
subgame-perfect equilibrium in this game. The ex-ante expected payoff of the
game is one half, but the first player has an advantage and gets greater
part of the pie. The uncertainty here arises because players initially do
not know who makes the first offer, and a natural measure of uncertainty is
the expectation of the squared difference between actual and expected
payoffs. 

The present paper generalizes this measure and illustrates with examples
that the uncertainty of the Shapley value is an interesting object that
behaves differently in different types of games. Some of its properties are
parallel to the properties of the Shapley value. The linearity property is,
however, replaced by convexity. It is also conjectured and proved for
certain games that the uncertainty of the sum of two games is larger than
the sum of the corresponding uncertainties.

The concept of uncertainty of a game outcome is not entirely new. %
\citeN{roth88}\footnote{%
See also \citeN{roth77b} and \citeN{roth77}.} has pointed out that the
Shapley value is insufficient for decision making purposes unless the
players are risk-neutral to a specific kind of uncertainty in the strategic
interactions, which he calls ``strategic risk''. He distinguishes this kind
of risk from ``ordinary risk'' that typically arise in lotteries, and he
characterizes the class of utility functions that may depend on this risk.
Roth, however, does not attempt to introduce a quantitative measure of
strategic risk, and that is the main contribution of the present paper.

The introduced measure provides a yardstick for both measuring the
reliability of the Shapley value and quantifying the strategic risk. The
main advantage of this measure is its simplicity. Future work will perhaps
provide alternative measures of uncertainty based on more accurate analysis
of bargaining randomness. The purpose of this paper is simply to bring
attention to the importance of measuring this randomness and to suggest a
possible solution.

The rest of the paper is organized as follows. Section 2 provides notation
and basic definitions. Section 3 lists examples. Section 4 shows certain
properties of the uncertainty concept, and Section 5 concludes.

\section{Notation and Definitions}

Let $A$ be a set of players: $A=\left\{ 1,2,...N\right\} .$ A game $G$ is a
map from a set of all subsets of $A$ to real numbers: $2^{A}\rightarrow 
\mathbb{R}$ . The interpretation of $G(S)$ is the payoff that the coalition
of players in subset $S$ can achieve on their own.

Let $(X,\Omega ,\mu )$ be a probability space and $f_{i}:X\rightarrow
2^{A-\{i\}}$ be a random variable that takes its values in the set of all
subsets of $A-\{i\},$ and has the following distribution function:%
\begin{equation}
\mu \left\{ f_{i}(X)=S\right\} =\frac{|S|!(N-|S|-1)!}{N!}.
\end{equation}%
Random variable $f_{i}$ can be interpreted as the random choice of coalition
that player $i$ joins. The particular choice of its probability distribution
is motivated by both the axiomatic approach and non-cooperative bargaining
games that support the Shapley value.

Define the marginal contribution of player $i$ as a function $d_{i}G$ that
acts on coalitions from $A-\{i\}$ in the following way:%
\begin{equation}
d_{i}G(S)=G(S+i)-G(S).
\end{equation}%
Then the Shapley value of game $G$ for player $i$ is the expectation of his
marginal contribution when the coalition is chosen randomly.%
\begin{equation}
V_{i}(G)=\mathbf{E}\left\{ d_{i}G\circ f_{i}\right\} .
\end{equation}%
The Shapley uncertainty of game $G$ for player $i$ is defined as the
variance of his marginal contribution:%
\begin{equation}
R_{i}(G)=\mathbf{Var}\left\{ d_{i}G\circ f_{i}\right\}
\end{equation}

This definition easily generalizes to other, non-Shapley, values. Indeed,
almost all reasonable values are probabilistic values (\citeN{weber88}),
which can be represented as the expectation of the player's marginal
contribution with respect to a certain probability distribution imposed on
coalitions. Of course, for different bargaining mechanisms the probability
distributions may be different. For example, another popular value, the
Banzhaf value, corresponds to the following probability distribution:%
\begin{equation}
\mu \left\{ f_{i}(X)=S\right\} =\frac{1}{2^{N-1}}.
\end{equation}

For probabilistic values, it is natural to define the corresponding
uncertainty as the variance of the marginal contribution. Most of the
properties that will be derived below for the uncertainty of the Shapley
value will also hold for the uncertainty of other probabilistic values. This
paper, however, sticks with the Shapley value as the most popular in
practical applications.

\section{Interpretation and Examples}

In many non-cooperative interpretations of the Shapley value, the actual
payoff of a player is a weighted average of his marginal contribution to a
random coalition and his expected payoff:%
\begin{equation}
\pi _{i}=\alpha (d_{i}G(S))+(1-\alpha )V_{i}(G).
\end{equation}%
For instance, in the Rubinstein-Gul-type bargaining models the parameter $%
\alpha $ is related to player impatience. Then, the variance of player
payoff is proportional to the Shapley uncertainty:%
\begin{equation}
\mathbf{Var}\left\{ \pi _{i}\right\} =\alpha ^{2}R_{i}(G).
\end{equation}

According to the Chebyshev inequality, it follows that 
\begin{equation}
\Pr \left\{ \left| \frac{\pi _{i}}{V_{i}(G)}-1\right| \geqslant c\right\}
\leqslant \alpha ^{2}\frac{R_{i}(G)}{\left[ V_{i}(G)\right] ^{2}},
\end{equation}%
that is, that the Shapley uncertainty provides a bound on the deviation of
the actual payoff from its expected value. Moreover, if the marginal
contributions are distributed approximately normally, which is likely to
happen in large games, then the statistical results for normal distributions
are applicable and the following inequality holds:%
\begin{equation}
\Pr \left\{ \left| \pi _{i}-V_{i}(G)\right| \geqslant 1.96\alpha \sqrt{%
R_{i}(G)}\right\} \leqslant 0.05
\end{equation}%
So, the Shapley uncertainty can be interpreted as measuring the variability
of the actual payoff.

Let me now present several examples.

\begin{example}
\label{example_additive_game}Additive game
\end{example}

In the additive game the value of a coalition is simply the sum of the
values of individual players. The Shapley value for a player coincides with
his individual value and, consequently, the risk is zero: Each player can
guarantee his own value.

\begin{example}
Uncertainty may be different for different players.
\end{example}

Let a three-player game be defined as follows. $G(\{1\})=0$, $G(\{2\})=3$, $%
G(\{3\})=6$, $G(\{1,2\})=24$, $G(\{2,3\})=18$, $G(\{1,3\})=21$, $%
G(\{1,2,3\})=60$. Then, 
\begin{equation*}
V_{1}=V_{2}=V_{3}=20,
\end{equation*}%
and 
\begin{equation*}
R_{1}=299;\;R_{2}=230;\;R_{3}=155.
\end{equation*}%
This example illustrates that players can have the same Shapley values but
different Shapley uncertainties.

\begin{example}
Majority game
\end{example}

In the majority game the value of a coalition is $N$ if the size of the
coalition is greater than half of the number of players, and $0$ otherwise.
For simplicity, let the number of players $N$ be odd. Then the Shapley value
of the game for player $i$ is 
\begin{equation}
V_{i}(G)=1,
\end{equation}%
and the Shapley uncertainty is%
\begin{equation}
R_{i}(G)=N.
\end{equation}%
So the uncertainty of the majority game increases with an increase in the
number of players. Intuitively this is because the marginal contribution of
a player is non-zero only if his vote is pivotal in a random ordering of
players, and in the majority game with a large number of players this
situation is very rare.

\begin{example}
\label{example_large_production}Large production economy \footnote{%
This example is based on an example of a production economy in %
\citeN{shapley_shubik67} and \citeN{osborne_rubinstein94}.}
\end{example}

An economy consists of $N$ agents, of which $kN$ are capitalists and ($1-k)N$
are workers. Each capitalist owns a factory and each worker owns one unit of
labor power. $K$ capitalists and $L$ workers can produce $\sqrt{KL}$ units
of output. What is the value and the uncertainty of the game for the
participants?

It is easy to check that for a random large coalition of size $X,$ the
frequency of capitalists is distributed according to the Gaussian law:%
\begin{equation}
\widehat{k}\sim \mathcal{N}\left( k,\frac{k}{X}\right) ,
\end{equation}%
where $\widehat{k}$ is the sample frequency of capitalists, and $\mathcal{N}%
\left( E,V\right) $ is the notation for the Gaussian distribution with
expectation $E$ and variance $V$.\ The marginal contribution of a worker
added to this coalition is 
\begin{equation}
d_{w}G(S)=\frac{1}{2}\sqrt{\frac{\widehat{k}}{1-\widehat{k}}}.
\end{equation}

For a large $X$ we can compute expectation and variance of this expression,
and averaging them over the coalition size, $X,$ gives the value and the
uncertainty of the game for a worker:%
\begin{eqnarray}
V_{w}(G) &=&\frac{1}{2}\sqrt{\frac{k}{1-k}}, \\
R_{w}(G) &\sim &\frac{1}{16(1-k)^{3}}\frac{\ln N}{N}\;(\text{where }%
N\rightarrow \infty ).
\end{eqnarray}

Similarly, for a capitalist we have 
\begin{eqnarray}
V_{c}(G) &=&\frac{1}{2}\sqrt{\frac{1-k}{k}}, \\
R_{c}(G) &\sim &\frac{1}{16k^{3}}\frac{\ln N}{N}\;(\text{where }N\rightarrow
\infty ).
\end{eqnarray}

In this simple example of production economy, the uncertainty of the game
decreases as the size of the game grows. The intuition behind this is
simple: As the number of participants grows, almost all coalitions have
approximately the same structure as the total population. For these typical
coalitions an agent has approximately the same marginal contribution, so his
uncertainty is low.

\begin{example}
Large Market Economy
\end{example}

Suppose there are just two commodities in the economy, apples and bread.
Assume there are $kN$ apple traders and $(1-k)N$ bread traders. Suppose also
that initially each trader has one unit of commodity of his own type, and
that the utility function of every trader is the same:%
\begin{equation}
u(a,b)=\sqrt{ab},
\end{equation}%
where $a$ is the number of apples and $b$ is the quantity of bread. Assume
that the utility is transferable.

Then a coalition of $A$ apple and $\ B$ bread traders has worth%
\begin{equation}
G(S)=\sqrt{AB},
\end{equation}%
and the problem is completely analogous to the problem of Example \ref%
{example_large_production}. Thus, for the apple trader the value and the
uncertainty are%
\begin{eqnarray}
V_{a}(G) &=&\frac{1}{2}\sqrt{\frac{1-k}{k}}, \\
R_{a}(G) &\sim &\frac{1}{16k^{3}}\frac{\ln N}{N}\;(\text{where }N\rightarrow
\infty ).
\end{eqnarray}

The expressions for the value and uncertainty of the bread trader are
similar. Therefore, this example suggests that in the market economy, as
well as in the production economy, the uncertainty goes to zero when the
size of the economy increases.

\begin{example}
Risk Factor Attribution
\end{example}

The Shapley values for the epidemiological example from %
\citeN{land_gefeller97} are $0.405,$ $0.08$, $0074$ respectively for
abnormal LDL, VLDL, and HDL cholesterol levels and $0.211$ for smoking. The
square roots of the corresponding Shapley uncertainties can be computed as $%
0.118,$ $0.056,$ $0.057,$ $0.11.$ The standard statistical test of
significance shows that only the effects of abnormal LDL and smoking are
significant at a 5 percent significance level.

\section{Properties of the Shapley Uncertainty}

Several of the properties of the Shapley uncertainty are similar to the
properties of the Shapley value and follow directly from the definition.
They are scaling, ``zero risk for dummy'', and symmetry properties:

\begin{proposition}[Scaling]
$R_{i}(tG)=t^{2}R_{i}(G)$
\end{proposition}

\begin{proposition}[Dummy]
If $G(X+i)=G(X)+G(i)$ for any $X$ that does not contain $i$ then $R_{i}(G)=0$%
.
\end{proposition}

\begin{proposition}[Symmetry]
If $G(X)=G(X-i+j)$ for any $X$ that contains $i$, then $R_{i}(G)=R_{j}(G)$.
\end{proposition}

For two-player games a stronger symmetry property holds. It states that in
two-player games both players are always equally exposed to the Shapley
uncertainty even if their values are different.

\begin{proposition}[Strong Symmetry]
\label{gamrisk:equality1} For any 2-player game $R_{i}(G)=R_{j}(G)$.
\end{proposition}

For games with more than two players this property is not true but still the
risks of different players cannot be too different, as the following
proposition shows.

\begin{proposition}
\label{gamrisk:inequality1} For every $N$-player game~$G$ 
\begin{equation*}
R_{i}(G)\leq (N-1)\sum_{j\not=i}R_{j}(G).
\end{equation*}
\end{proposition}

\textbf{\noindent Proof:} First, the Shapley uncertainty can be defined in
an equivalent way. Let $\mathfrak{B}$ be the space of all orderings of
elements in $A.$ Define a map $p_{i}:\mathfrak{B\rightarrow }2^{A}$ that
takes an ordering to the subset of elements of $A$ that precedes element $i$
in the ordering, and let $\widetilde{d_{i}G}=$ $d_{i}G\circ p_{i}.$Take a
random variable $g:X\rightarrow \mathfrak{B}$ that has a uniform
distribution on orderings. Then it is easy to check that: 
\begin{eqnarray}
f_{i} &=&p_{i}\circ g, \\
V_{i}(G) &=&\mathrm{E}(\widetilde{d_{i}G}\circ g),  \notag \\
R_{i}(G) &=&\limfunc{Var}(\widetilde{d_{i}G}\circ g).  \notag
\end{eqnarray}%
Second, for any ordering $b\in \mathfrak{B}$, the sum of marginal
contributions of all players is equal to the value of the game, which is not
random:%
\begin{equation}
\sum_{i}\widetilde{d_{i}G}(b)=G(A).
\end{equation}%
It follows 
\begin{equation}
\mathrm{Var}(\widetilde{d_{i}G}\circ g)=\mathrm{Var}(\sum_{j\not=i}%
\widetilde{d_{j}G}\circ g)\leq (N-1)\sum_{i\not=j}\mathrm{Var}(\widetilde{%
d_{j}G}\circ g),
\end{equation}%
where the last inequality is a well-known inequality for the variance of the
sum of random variables. So, 
\begin{equation}
R_{i}(G)\leq (N-1)\sum_{j\not=i}R_{j}(G).
\end{equation}%
QED

This result suggests that the uncertainty of a game has a component that is
independent of the strategic situation of an individual. In large games the
results of similar type imply that even if a player does not know his exact
role in a game, he is able to evaluate the probable range of the risk that
he is going to encounter. Therefore, the player could meaningfully compare
games on the basis of their uncertainty, or riskiness. A possible measure of
expected uncertainty of a game is the average of the players' uncertainties:%
\begin{equation}
R(G)=\frac{1}{N}\sum_{i=1}^{N}R_{i}(G).
\end{equation}
It is the uncertainty that a player can expect to encounter if he is not yet
aware of his role in the game.

The properties above have been mostly concerned with comparison of the
Shapley uncertainties for the different players. Another important property,
convexity, is concerned with behavior of the uncertainty relative to the
addition of games.

\begin{proposition}[Convexity]
\label{gamrisk:risk_upper_bound} For each pair of games $G$ and $H$, 
\begin{equation*}
R_{i}(\alpha G+(1-\alpha )H)\leq \alpha R_{i}(G)+(1-\alpha )R_{i}(H).
\end{equation*}%
The equality is only possible if the games are proportionate to each other.
\end{proposition}

\textbf{\noindent Proof:} The uncertainty of a game is the variance of a
random variable. The assertion of the proposition follows from the
well-known property of the variance of a weighted sum of two random
variables. QED

From this property the following upper bound on the Shapley uncertainty of
the sum of two games follows.

\begin{corollary}
\bigskip \label{gamrisk:risk_upper_bound2} For each pair of games $G$ and $H$%
, 
\begin{equation*}
R_{i}(G+H)\leq 2\left[ R_{i}(G)+R_{i}(H)\right] .
\end{equation*}%
The equality is only possible if the games are proportionate to each other.
\end{corollary}

Corollary \ref{gamrisk:risk_upper_bound2} gives an upper bound to the
uncertainty of the sum of two games. What can be said about the lower bound?
The question boils down to whether two games with large Shapley uncertainty
can add up to a game with low uncertainty. Numerical experiments suggest the
following conjecture:

\begin{conjecture}
There exist such $d(n)>0$ that for each pair of superadditive $n$-player
games $G$ and $H$ 
\begin{equation*}
R_{i}(G+H)\geq d(n)(R_{i}(G)+R_{i}(H)).
\end{equation*}
\end{conjecture}

This conjecture is true for symmetric convex games. In these games, the
worth of a coalition depends only on its size and the dependence is
non-decreasing and convex:%
\begin{equation}
G(S)=g(|S|),\text{ where }g^{\prime }(x)\geq 0\text{, and }g^{\prime \prime
}(x)\geq 0
\end{equation}%
For these games, the following proposition holds

\begin{proposition}
For each pair of symmetric convex games $G$ and $H$,%
\begin{equation*}
R_{i}(G+H)\geq R_{i}(G)+R_{i}(H).
\end{equation*}
\end{proposition}

\textbf{\noindent Proof: }Clearly the proposition holds if player $i$'s
marginal contributions have positive covariance for each pair of symmetric
convex games $G$ and $H:$%
\begin{equation}
\mathrm{Cov}\{d_{i}G(S),d_{i}H(S)\}\geq 0,
\end{equation}%
where $S$ is the random coalition. In other words, it is sufficient to prove
that 
\begin{equation}
\mathrm{E}\{d_{i}G(S)d_{i}H(S)\}-\mathrm{E}\{d_{i}G(S)\}\mathrm{E}%
\{d_{i}H(S)\}\geq 0.  \label{positive_covariance}
\end{equation}%
Marginal contributions $d_{i}G(S)$ and $d_{i}H(S)$ are non-negative and
non-decreasing in argument $S,$ and therefore inequality (\ref%
{positive_covariance}) holds because of the following lemma:

\begin{lemma}
If $f_{1}(x)$ and $f_{2}(x)$ are two non-negative, non-decreasing functions,
and $\mathcal{F}(x)$ is the cumulative distribution function of a random
variable, then 
\begin{equation}
\int_{-\infty }^{\infty }f_{1}(x)f_{2}(x)d\mathcal{F}(x)\geq \left(
\int_{-\infty }^{\infty }f_{1}(x)d\mathcal{F}(x)\right) \left( \int_{-\infty
}^{\infty }f_{2}(x)d\mathcal{F}(x)\right) .
\end{equation}
\end{lemma}

\textbf{\noindent Proof of Lemma: }Note that if the lemma holds for $%
\{f_{1}^{\prime },f_{2}\}$ and $\{f_{1}^{\prime \prime },f_{2}\}$, then it
also holds for $\{\alpha ^{\prime }f_{1}^{\prime }+\alpha ^{\prime \prime
}f_{1}^{\prime \prime },f_{2}\},$ where $\alpha ^{\prime }$ and $\alpha
^{\prime \prime }$ are non-negative coefficients. Similarly, if the lemma
holds for $\{f_{1},f_{2}^{\prime }\}$ and $\{f_{1},f_{2}^{\prime \prime }\}$%
, then it also holds for $\{f_{1},\beta ^{\prime }f_{2}^{\prime }+\beta
^{\prime \prime }f_{2}^{\prime \prime }\},$ where $\beta ^{\prime }$ and $%
\beta ^{\prime \prime }$ are non-negative. It follows that the lemma needs
only to be proved for a set of elementary functions that can approximate all
non-negative, non-decreasing functions by linear combinations with
non-negative coefficients. It is enough to take the set of the following
elementary functions $\phi _{T}:$%
\begin{equation}
\phi _{T}(x)=\left\{ 
\begin{array}{cc}
1 & \text{if }x\geq T, \\ 
0 & \text{if }x<T,%
\end{array}%
\right. 
\end{equation}%
which are characteristic functions of sets $[T;\infty ).$

Define measure $\mu _{\mathcal{F}}$ as follows:%
\begin{equation}
\mu _{\mathcal{F}}\{X\}=\int_{-\infty }^{\infty }\chi _{X}(x)d\mathcal{F}(x),
\end{equation}%
where $\chi _{X}(x)$ denotes the characteristic function of set $X.$

In terms of this measure, the inequality for elementary functions becomes:%
\begin{equation}
\mu _{\mathcal{F}}\{[T_{1},\infty )\cap \lbrack T_{2},\infty )\}\geq \mu _{%
\mathcal{F}}\{[T_{1},\infty )\}\mu _{\mathcal{F}}\{[T_{2},\infty )\},
\end{equation}%
which is evidently true because $[T_{1},\infty )\cap \lbrack T_{2},\infty )$
is either $[T_{1},\infty )$ or $[T_{2},\infty ),$ and $\mu _{\mathcal{F}%
}\{X\}\leq 1$ for any $X.$ QED.

\section{Conclusion}

In the course of this paper, the concept of the Shapley uncertainty has been
rigorously defined and illustrated by means of several examples. Certain
properties of this concept have been proved that show how the uncertainty
changes under scaling in the size of the game, substitution in the role of
players, and taking an average of games. These examples and propositions
illustrate that the Shapley uncertainty is an object, that has many
interesting properties. The main question is, however: What is this concept
good for?

Actually, it is good for several important things. First, it measures the
reliability of Shapley value solutions to various practical problems.
Second, it is useful in checking the robustness of theoretical arguments
that employ the Shapley value as the embodiment of the outcome in market
games. Third, the Shapley uncertainty is an additional dimension that a
player should take into account if he evaluates the prospect of playing a
game. Finally, the Shapley uncertainty may be helpful in the design of games
as implementation mechanisms, where the designer is interested in ensuring
the stability of the outcome and making the mechanism attractive to
risk-sensitive participants.

\bibliographystyle{CHICAGO}
\bibliography{comtest}

\end{document}